\documentclass[11pt,reqno]{amsart}
\usepackage[T1]{fontenc}
\usepackage{lmodern,amssymb,microtype}
\usepackage[letterpaper,margin=1in]{geometry}
\numberwithin{equation}{section}
\usepackage[hidelinks]{hyperref}
\hypersetup{pdftitle={A Finiteness Result for Commuting Squares and Quantitative Ocneanu Rigidity},
  pdfauthor={Shannon Hall and Remus Nicoar\u{a}}}
\newtheorem{theorem}{Theorem}[section]
\newtheorem{proposition}[theorem]{Proposition}
\newtheorem{corollary}[theorem]{Corollary}

\theoremstyle{definition}
\newtheorem{definition}[theorem]{Definition}
\theoremstyle{remark}
\newtheorem{remark}[theorem]{Remark}
\newcommand{\C}{\mathbb C}
\newcommand{\U}{\mathcal U}
\newcommand{\Tr}{\operatorname{Tr}}
\newcommand{\spn}{\operatorname{span}_{\C}}
\newcommand{\End}{\operatorname{End}}
\title[Finiteness for Commuting Squares and Quantitative Ocneanu Rigidity]{A Finiteness Result for Commuting Squares and Quantitative Ocneanu Rigidity}
\author{Shannon Hall}
\author{Remus Nicoar\u a}
\makeatletter
\g@addto@macro\@setauthors{%
  {\centering\footnotesize
   Department of Mathematics, University of Tennessee, Knoxville\\
   1403 Circle Drive, Knoxville, TN 37996, USA\\[3pt]
   Corresponding author: Remus Nicoar\u a,
   \href{mailto:rnicoara@utk.edu}{\texttt{rnicoara@utk.edu}}\par}}
\makeatother
\date{November 26, 2026}
\begin{document}
\begin{abstract}
We give a commuting-square proof of Ocneanu's finiteness theorem for standard invariants of subfactors with a prescribed finite principal graph. The proof gives an explicit upper bound on the number of isomorphism classes in terms of the graph alone. We introduce the spanning relative commutant condition for nondegenerate Markov commuting squares of finite-dimensional algebras and show that it characterizes depth at most two for the column subfactor of a connected commuting square. For a fixed upper-right algebra of complex dimension $D\geq 2$ and a fixed trace, there are at most $(3D)^{25D}$ isomorphism classes of commuting squares satisfying this condition.
\end{abstract}
\maketitle

\section{Introduction}

Popa introduced the commuting-square condition in his 1983 paper on maximal
injective subalgebras of free group factors \cite[Lemma~1.2.2]{MaxInj}.
His paper on orthogonal pairs of subalgebras, published the same year, treats
the case of scalar intersection \cite[Lemma~2.1 and Definition~2.2]{Orth}. The terminology ``commuting square''
was introduced later by Jones, as Pimsner and Popa note \cite[\S1.4, p.~279]{PPApprox}.
It appears in the book of Goodman, de la Harpe, and Jones \cite[Definition~4.1.1]{GHJ}.
By a matrix algebra we mean a finite direct sum of full matrix algebras
over $\C$. Let $P_0$ be such an algebra with a faithful normalized trace
$\tau$, and let $P_{-1}$, $Q_0$ and $Q_{-1}$ be unital $*$-subalgebras
of $P_0$ with $Q_{-1}\subset P_{-1}\cap Q_0$. Consider the diagram
\begin{equation}\label{eq:square}
 \mathcal C=
 \left(\begin{array}{ccc}
  P_{-1}&\subset&P_0\\
  \cup&&\cup\\
  Q_{-1}&\subset&Q_0
 \end{array},\,\tau\right).
\end{equation}
For a unital $*$-subalgebra $X\subset P_0$, let $E_X$ denote the
$\tau$-preserving conditional expectation onto $X$. We call
\eqref{eq:square} a \emph{commuting square} if
\[
 E_{P_{-1}}E_{Q_0}=E_{Q_0}E_{P_{-1}}=E_{Q_{-1}} \qquad\text{on }P_0.
\]
Equivalently, $P_{-1}\ominus Q_{-1}$ and $Q_0\ominus Q_{-1}$ are
orthogonal in $L^2(P_0,\tau)$.

Following Popa \cite[\S1.1.5, p.~172]{Pa}, we call the commuting square
$\mathcal C$ \emph{nondegenerate} if
\[
 P_0=P_{-1}\cdot Q_0
    =\operatorname{span}_{\C}\{xy:x\in P_{-1},\ y\in Q_0\}.
\]

We assume throughout that all inclusions are unital, that this condition holds, and that $\tau$ is a
$\lambda$-Markov trace for both horizontal inclusions, with the same scalar
$\lambda>0$.

Suppose, in addition, that all four inclusion graphs are connected and that
the Markov traces are compatible in both directions.
If the common vertical Markov index is greater than one, iterating the
basic construction along the columns $Q_{-1}\subset P_{-1}$ and
$Q_0\subset P_0$ gives two towers whose weak closures are hyperfinite $\mathrm{II}_1$
factors $N$ and $M$, respectively. We call $N\subset M$ the
\emph{column subfactor} of $\mathcal C$
\cite[Ch.~5]{JS}, \cite[Remark~A.15]{BCEGP}. 

For a finite-index subfactor $N=M_{-1}\subset M=M_0$, let
$(M_j)_{j\geq-1}$ denote its Jones tower \cite{J83}. Its principal graph
is the rooted bipartite multigraph associated to the tower
\[
  \C=N'\cap N\ \subset\ N'\cap M\ \subset\ N'\cap M_1\ \subset\ \cdots .
\]
The subfactor has \emph{finite depth} if its principal graph is finite,
a notion introduced by Ocneanu in his work on paragroups in the mid-1980s
\cite{Oc,Sa}.
Its \emph{depth} is then the greatest distance from the root to a vertex
of the graph. The depth of a commuting square is the depth of its column
subfactor.

Applying Jones's basic construction \cite{J83} to the two rows of \eqref{eq:square} gives their common extension \cite[\S1.1.6]{Pa}, which we denote by $\mathcal L$ as in \cite{Ni}. Thus
\begin{equation}\label{grid}
 \mathcal L=
 \left(\begin{array}{ccccc}
  P_{-1}&\subset&P_0&\overset{e}{\subset}&P_1\\
  \cup&&\cup&&\cup\\
  Q_{-1}&\subset&Q_0&\overset{e}{\subset}&Q_1
 \end{array},\,\tau\right),
 \qquad P_1=\langle P_0,e\rangle,\quad Q_1=\langle Q_0,e\rangle .
\end{equation}
Here $e$ is the common Jones projection, and $\tau$ also denotes the extended trace on $P_1$. With our convention, the Markov property reads $\tau(xey)=\lambda\tau(xy)$ for $x,y\in P_0$, or equivalently $E_{P_0}(e)=\lambda 1$. We write $\U(X)$ for the unitary group of $X$ and $\|x\|_2=\tau(x^*x)^{1/2}$.

\begin{definition}\label{def:src}
Put $R_0=P_{-1}'\cap Q_0$ and $R_1=P_{-1}'\cap Q_1$ in \eqref{grid},
with commutation taken in $P_1$.
We say that $\mathcal C$ has a \emph{spanning relative commutant} \textup{(SRC)} if
\begin{equation*}\tag{SRC}\label{eq:src}
                         R_1\cdot Q_0=Q_1 .
\end{equation*}
\end{definition}

Since $e\in R_1$ and $Q_1=\langle Q_0,e\rangle$, the algebras $R_1$ and $Q_0$ generate $Q_1$ for every nondegenerate commuting square. The basic construction writes each element of $Q_1$ as a sum of terms $qeq'$, with coefficients from $Q_0$ on both sides of $e$. Condition \eqref{eq:src} requires an expansion
\[
               x=\sum_j r_jq_j,\qquad r_j\in R_1,\quad q_j\in Q_0
               \qquad(x\in Q_1).
\]
Thus $R_1$ spans $Q_1$ as a right $Q_0$-module.

\begin{theorem}[Depth $2$ criterion]\label{th:src-depth}
Suppose that $\mathcal C$ gives a subfactor $N\subset M$ by the column
construction above. Then
\[
       R_1\cdot Q_0=Q_1
       \quad\Longleftrightarrow\quad
       \text{the principal graph of }N\subset M\text{ has depth at most two}.
\]
For a proper irreducible inclusion this is equivalent to depth exactly two.
\end{theorem}

The proof in Section~\ref{sec:expectation} identifies the orthogonal
projection $z$ onto $R_1\cdot Q_0$ with the unit of the basic-construction
ideal in $N'\cap M_2$. In the decomposition of $N'\cap M_2$ into full
matrix algebras, $1-z$ is the sum of the units of the summands indexed
by vertices at distance three from the root of the principal graph.
The dimension formula \eqref{eq:size} gives
another test for SRC. In particular, when $R_0=\C1$ it becomes
$\dim_\C Q_1=(\dim_\C R_1)(\dim_\C Q_0)$.

An isomorphism of two commuting squares is a trace-preserving $*$-isomorphism between their upper-right algebras which maps each of the other three algebras onto the corresponding algebra. Its unitary extension to the trace Hilbert spaces intertwines the Jones projections, hence carries $Q_1$ and $R_1$ onto their counterparts. In particular SRC is invariant under isomorphism.

Fix the four corner algebras, the trace on $P_0$, and the inclusion
matrix of each edge, using consistent labels for the simple summands
at each corner. For these fixed \emph{tracial inclusion data}, we prove
that commuting squares satisfying SRC are isolated up to isomorphism
and that only finitely many isomorphism classes occur. Finiteness
also holds when only $(P_0,\tau)$ is fixed.

\begin{theorem}\label{th:fin}
\textup{(i)} Let $\mathcal C$ be a nondegenerate commuting square with Markov trace $\tau$ as in \eqref{eq:square}--\eqref{grid}, and let $(u_m)_{m\geq1}$ be a sequence in $\U(Q_{-1}'\cap P_0)$ with $\|u_m-1\|_2\to0$. For each $m\geq1$, suppose that
\begin{equation}\label{eq:moving-square}
 \mathcal C_m=
 \left(\begin{array}{ccc}
      P_{-1}&\subset&P_0\\
      \cup&&\cup\\
      Q_{-1}&\subset&Q_0^{(m)}=u_mQ_0u_m^*
 \end{array},\,\tau\right)
\end{equation}
is a commuting square. Put $Q_1^{(m)}=\langle Q_0^{(m)},e\rangle\subset P_1$ and $R_1^{(m)}=P_{-1}'\cap Q_1^{(m)}$. If
\[
                   R_1^{(m)}\cdot Q_0^{(m)}=Q_1^{(m)}
                          \qquad\text{for every }m\geq1,
\]
then there is an $m_0$ such that $\mathcal C_m$ is isomorphic to $\mathcal C$ for every $m\geq m_0$.

\smallskip
\textup{(ii)} For fixed tracial inclusion data, there are only finitely many isomorphism classes of nondegenerate commuting squares satisfying \eqref{eq:src}. The same finiteness conclusion holds if only $(P_0,\tau)$ is fixed. For each integer $D\geq1$ it also holds for all such commuting squares with $\dim_\C P_0\leq D$ and $Z(P_{-1})\cap Z(P_0)=\C1$.
\end{theorem}

For a general sequence with fixed data, identify $Q_{-1}\subset P_{-1}\subset P_0$ and take convergence of the fourth corners to mean $E_{Q_0^{(m)}}\to E_{Q_0}$ in operator norm on $L^2(P_0,\tau)$. Section~\ref{sec:finiteness} shows how to put a subsequence in the form \eqref{eq:moving-square} with $u_m\to1$.

To prove \textup{(i)}, suppose that infinitely many $\mathcal C_m$ are not isomorphic to $\mathcal C$. After normalizing the corresponding unitaries and passing to a subsequence, we obtain the self-adjoint unit vector
\[
                    h=\lim_{m\to\infty}
                       \frac{u_m-1}{i\|u_m-1\|_2},
                  \qquad h=h^*,\qquad \|h\|_2=1 .
\]
As in \cite[Lemma~1.8]{N06} and \cite[Proposition~2.9]{Ni}, this normalization makes $h$ orthogonal to each of the four algebras in \eqref{eq:direction-orthogonal}. Applying Proposition~\ref{pr:expectation} to each $\mathcal C_m$ gives a conditional expectation $F_m:P_1\to P_0$ with
\[
             F_m(Q_1^{(m)})=Q_0^{(m)},\qquad
             F_m\big((R_1^{(m)})'\cap P_1\big)=P_{-1}.
\]
Passing to a subsequence of the pairs $(R_1^{(m)},F_m)$ gives the two inclusions in \eqref{eq:limitimages}. Equation \eqref{eq:direction-sum} then puts $h$ in the sum of the four algebras in \eqref{eq:direction-orthogonal}, forcing $h=0$. Compactness of the unitary parameter space then gives \textup{(ii)}.

The \emph{standard invariant} of $N\subset M$ is the lattice of higher
relative commutants $M_i'\cap M_j$, $-1\leq i\leq j$, with its inclusions,
traces and distinguished Jones projections, up to isomorphism preserving
these data \cite[\S1.3]{Pa}, \cite[\S1]{PopaAxioms}.

For $N\subset M$ of finite depth and all sufficiently large $k$, the diagram
\[
 \begin{array}{ccc}
 M'\cap M_k   &\subset& M'\cap M_{k+1}\\
 \cup         &       & \cup\\
 M_1'\cap M_k &\subset& M_1'\cap M_{k+1}
 \end{array}
\]
with the trace inherited from $M_{k+1}$ is a \emph{standard commuting square}
for $N\subset M$. Both rows continue from these columns by Jones's basic
construction. By Popa's reconstruction theorem, the trace-preserving isomorphism
class of this commuting square determines the entire standard invariant of $N\subset M$
(see \cite{Po} and \cite[\S1.3.7, p.~185]{Pa}).

For a subfactor satisfying the depth-two condition, Proposition~\ref{pr:depth}
verifies SRC directly for the standard commuting square \eqref{eq:standard-square}.
This was our original motivation for the condition. It extends the
finiteness result for large second relative commutants
\cite[Corollary~5.2]{Ni} to these standard commuting squares. Applied after grouping
levels in the Jones tower, Theorem~\ref{th:fin} gives the following
finiteness theorem of Ocneanu.

\begin{theorem}[Ocneanu]\label{th:oc}
Let $\Gamma$ be a finite rooted bipartite multigraph. Among finite-index inclusions $N\subset M$ of $\mathrm{II}_1$ factors whose principal graph is isomorphic to $\Gamma$ as a rooted graph, only finitely many isomorphism classes of standard invariants occur.
\end{theorem}

Ocneanu gave his original proof in the language of paragroups, which he introduced in \cite{Oc}. As Kawahigashi explains in \cite[Theorem~2.1 and Corollary~2.2]{Kaw}, Ocneanu first proved rigidity for fixed fusion rules and then deduced finiteness for a prescribed principal graph. In Section~\ref{sec:depth} we prove the principal-graph statement by passing from $N\subset M$ to $N\subset M_{\ell-1}$, where $\ell$ is the greatest distance from the root to an odd vertex. We show that this inclusion satisfies the depth~$2$ condition of Section~\ref{sec:depth} and hence that its standard commuting square satisfies \eqref{eq:src}. We recover the original standard invariant by taking a corner. Corollary~\ref{cor:oc-rigidity} then deduces the fixed-fusion-rule statement from Theorem~\ref{th:oc}.

\subsection*{Quantitative estimates}
Section~\ref{sec:quantitative} gives a common bound for the two finiteness theorems. Put $B(1)=1$ and $B(D)=(3D)^{25D}$ for $D\ge2$. By Theorem~\ref{th:quant-squares}, there are at most $B(D)$ isomorphism classes of SRC commuting squares when $\dim_\C P_0\le D$ and either $(P_0,\tau)$ is prescribed or the upper horizontal inclusion is connected, that is, $Z(P_{-1})\cap Z(P_0)=\C1$. The proof expresses SRC as a rank condition on a matrix depending on the unitary that parametrizes the commuting square. The isomorphism class is locally constant on each set of solutions of fixed rank. Milnor's estimate \cite[Theorem~2]{Mil} bounds the number of its connected components.

Let $\Gamma$ be a prescribed finite principal graph with adjacency matrix $G$ and root $*$. Let $\ell$ be the greatest distance from $*$ to an odd vertex and put $D_\Gamma=(G^{6\ell})_{*,*}$. Theorem~\ref{th:quant-graph} bounds the number of isomorphism classes of standard invariants with principal graph $\Gamma$ by $B(D_\Gamma)$. To apply the bound for commuting squares, we use the standard commuting square \eqref{eq:grouped-square} of $N\subset M_{\ell-1}$. Its upper-right algebra is $N'\cap M_{3\ell-1}$, whose dimension is $D_\Gamma$, the number of closed walks of length $6\ell$ from the root. For each such commuting square, we then count the possible block ranks of the projection in $N'\cap M_{\ell-1}$ that recovers the standard invariant of $N\subset M$. Section~\ref{sec:quantitative} also gives a bound when only the index and depth are bounded.

For spin models of order greater than one, Theorems~\ref{th:src-depth}
and~\ref{th:spin} identify the matrices giving depth two as the character
tables of finite abelian groups, up to Hadamard equivalence. Proposition~\ref{pr:spin-separation} gives a separation estimate for these commuting squares. If $A$ and $B$ are the fourth corners of two such commuting squares of order $n$ satisfying SRC, with common upper-left algebra $D_n$, then
\[
              \sup_{\|x\|_2=1}\|E_A(x)-E_B(x)\|_2<\frac1n
       \quad\Longrightarrow\quad
              \text{the two commuting squares are isomorphic}.
\]
Here $x\in M_n(\C)$ and $\|x\|_2=\tau(x^*x)^{1/2}$, where $\tau$ is the normalized matrix trace.

In developing the proofs, we used OpenAI's ChatGPT to explore possible proof strategies and to review arguments and references.

\section{The depth 2 criterion and a conditional expectation}\label{sec:expectation}

\subsection*{The relative commutants of the column subfactor}
Write $Q_i\subset P_i$ for column $i$ of the horizontal extension of
\eqref{grid}. Put $A_{0,i}=Q_i$ and $A_{1,i}=P_i$, and continue each
column by the basic construction to obtain $A_{k,i}$ for $k\geq 2$.
Taking weak closures in the trace representation gives the Jones tower
\[
 M_i=\Big(\bigcup_{k\geq0} A_{k,i}\Big)'',\qquad
 N=M_{-1}\subset M=M_0\subset M_1\subset M_2\subset\cdots .
\]
The common horizontal Jones projections $e_i\in Q_i$ are the Jones
projections in this tower, with $e_1=e$ as in \eqref{grid}. In particular
$[M:N]=\tau(e_1)^{-1}=\lambda^{-1}$. For
$i\geq0$ put $R_i=P_{-1}'\cap Q_i$ inside $P_i$. Ocneanu compactness,
applied to the rectangle from column $-1$ to column $i$, gives
\begin{equation}\label{eq:compact-columns}
                    R_i=N'\cap M_i .
\end{equation}
We use the general form of compactness in \cite[Theorem~C.1]{BCEGP}. See also \cite[Theorem~5.7.1]{JS}.
Thus $R_0,R_1,R_2$ are the first three higher relative commutants
of the associated subfactor.

\begin{proof}[Proof of Theorem~\ref{th:src-depth}]
Represent $Q_2=\langle Q_1,e_2\rangle$ on $\mathcal H=L^2(Q_1)$.
By \eqref{eq:compact-columns}, $I=\spn R_1e_2R_1$ is the
basic-construction ideal in $R_2$. Its unit $z$ is the orthogonal
projection onto $I\mathcal H$. Since $1\in R_1$ and $e_2$ projects
onto $L^2(Q_0)$, we have
\[
                I\mathcal H=R_1\cdot(e_2\mathcal H)=R_1\cdot Q_0 .
\]
Thus SRC holds precisely when $z=1$, or equivalently $I=R_2$.
The full matrix summands of $(1-z)R_2$ are indexed by the vertices
at distance three from the root \cite[\S1.3]{Pa}. Their absence is equivalent
to depth at most two, since every path from the root to a deeper vertex
passes through depth three.
For a proper irreducible inclusion a principal graph of depth one would
consist of a single edge and have index one. This proves the last assertion.
\end{proof}

More generally, fix $j\geq1$ and represent $Q_{j+1}$ on $L^2(Q_j)$.
By \eqref{eq:compact-columns}, $I_j=\spn R_j e_{j+1}R_j$ is the
basic-construction ideal in $R_{j+1}$. Since $1\in R_j$ and $e_{j+1}$
projects onto $L^2(Q_{j-1})$, the range of the unit of $I_j$ is
\[
 I_j L^2(Q_j)=R_j\cdot Q_{j-1}.
\]
Thus $R_j\cdot Q_{j-1}=Q_j$ exactly when $I_j=R_{j+1}$. The full matrix summands of
$R_{j+1}$ outside $I_j$ correspond to the principal graph vertices at
distance $j+2$ from the root \cite[\S1.3]{Pa}. Hence
\begin{equation}\label{eq:higher-src}
 (\mathrm{SRC}_j)\qquad
 R_j\cdot Q_{j-1}=Q_j
 \quad\Longleftrightarrow\quad
 \operatorname{depth}(N\subset M)\leq j+1 .
\end{equation}
In particular the subfactor has finite
depth exactly when one of these equalities holds.

\subsection*{The dimension test and the expectation}
We return to the general finite-dimensional setting of \eqref{grid}.
The dimension of $R_1\cdot Q_0$ can be read from the inclusions
$R_0\subset R_1$ and $R_0\subset Q_0$.

Since $Q_1$ is spanned by $Q_0eQ_0$ and $E_{P_0}(e)=\lambda 1$, we have $E_{P_0}(Q_1)=Q_0$. Also $E_{P_0}(R_1)\subset R_0$ by $P_{-1}$-bimodularity, so that $E_{Q_0}|_{R_1}=E_{R_0}|_{R_1}$. Choose matrix units $(f^{(i)}_{ab})$ in each simple block of $R_0$ and put $t_i=f^{(i)}_{11}$. Choose a basis $(v_{i\alpha})_\alpha$ of $R_1t_i$ orthonormal for the inner product $\tau(y^*x)/\tau(t_i)$. Then
\[
 E_{Q_0}(v_{i\alpha}^*v_{j\beta})
   =\delta_{ij}\delta_{\alpha\beta}\,t_i,
\]
since $t_iR_0t_j=0$ for $i\ne j$ and $t_iR_0t_i=\C t_i$. Consequently left multiplication by $v_{i\alpha}$ maps $t_iQ_0$ isometrically onto $v_{i\alpha}t_iQ_0$ in $L^2(Q_1,\tau)$, and these ranges are mutually orthogonal. They span $R_1\cdot Q_0$ because
\[
 rq=\sum_{i,a}(r f^{(i)}_{a1})(f^{(i)}_{1a}q)
       \qquad(r\in R_1,\ q\in Q_0).
\]
Counting dimensions gives
\begin{equation}\label{eq:size}
       R_1\cdot Q_0=Q_1\quad\Longleftrightarrow\quad
       \dim_\C Q_1=\sum_i\dim_\C(R_1t_i)\,\dim_\C(t_iQ_0).
\end{equation}
When $R_0=\C1$, the commuting square satisfies \eqref{eq:src} if and only if $\dim_\C Q_1=(\dim_\C R_1)(\dim_\C Q_0)$.

The isolation proof also needs an expectation onto $P_0$ which sends
the commutant of $R_1$ into $P_{-1}$. We construct it by weighting
the trace-preserving expectation.

\begin{proposition}\label{pr:expectation}
Let $\mathcal L$ be the basic construction \eqref{grid} of the commuting square $\mathcal C$ in \eqref{eq:square}, and put $R_1=P_{-1}'\cap Q_1$. Suppose that $R_1\cdot Q_0=Q_1$. There is a positive invertible $d\in P_0'\cap Q_1$ with $E_{P_0}(d)=1$ such that $F(x)=E_{P_0}(dx)$, $x\in P_1$, defines a faithful conditional expectation of $P_1$ onto $P_0$ satisfying
\begin{equation}\label{eq:expectation-images}
                         F(Q_1)=Q_0,\qquad F(R_1'\cap P_1)=P_{-1}.
\end{equation}
\end{proposition}

\begin{proof}
Represent $P_0$ by left multiplication on $\mathcal H=L^2(P_0,\tau)$, and denote right multiplication by $\rho$. In this proof absolute commutants are taken in $\mathcal B(\mathcal H)$, so $P_1=\rho(P_{-1})'$. Since $P_0=P_{-1}\cdot Q_0$ and $e$ commutes with $P_{-1}$, the identity $P_1=\spn P_0eP_0$ gives $P_1=P_0\cdot Q_1$. Taking adjoints in \eqref{eq:src} gives $Q_1=Q_0\cdot R_1$. It follows that
\[
                         P_1=P_0\cdot R_1=R_1\cdot P_0.
\]
Together with $E_{P_0}(R_1)\subset R_0$, this makes the inclusions
$R_0\subset P_0$ and $R_1\subset P_1$ a nondegenerate commuting square.

Let $J\xi=\xi^*$ on $L^2(P_0,\tau)$. We claim that $JR_1J=R_1$. Indeed, $R_1$ consists of the operators commuting with both actions of $P_{-1}$ and carrying $Q_0$ into $Q_0$. Every element of $Q_1=\langle Q_0,e\rangle$ preserves $Q_0$ and commutes with the right action of $P_{-1}$. Conversely, the restriction to $Q_0$ of an operator with the stated properties is right $Q_{-1}$-linear. It thus belongs to $\End_{Q_{-1}}(Q_0)$, the restriction of $Q_1$ to $L^2(Q_0,\tau)$ by the basic construction. Right $P_{-1}$-linearity determines the operator on $Q_0\cdot P_{-1}=P_0$. This description of $R_1$ is invariant under $J$.

Conjugating this commuting square by $J$ gives
\begin{equation}\label{eq:reflected-square}
 \left(\begin{array}{ccc}
   R_1&\subset&L=P_{-1}'\\
   \cup&&\cup\\
   I=\rho(R_0)&\subset&K=\rho(P_0)
 \end{array},\,\tau^J\right),
 \qquad L=K\cdot R_1 .
\end{equation}
Expectations in \eqref{eq:reflected-square} are taken with respect to $\tau^J(y)=\overline{\tau(JyJ)}$ for $y\in L$. In particular $I\subset R_1\cap P_0'$. Choose orthonormal bases $(f_l)$ of $I$ and $(u_k)$ of $K$ for $\tau^J$.
Set $a=(\sum_l f_l^*f_l)^{-1}$. The sum is independent of the basis.
Using normalized matrix units shows that on a block $M_n$ of $I$ with
minimal trace $\gamma$, it equals $(n/\gamma)1$. In particular $a$ is
positive and invertible in $Z(I)$. Put $v_k=u_ka^{1/2}$. Expanding
$E_I$ in $(f_l)$ and using Parseval's identity in $K$ gives the right
Pimsner--Popa basis identity
\[
   \sum_k v_k E_I(v_k^*y)
    =\sum_l y f_l^* a f_l=y\qquad(y\in K).
\]

Since $E_{R_1}|_K=E_I|_K$ in \eqref{eq:reflected-square} and $L=K\cdot R_1$, the same $(v_k)$ is a right basis of $L$ over $R_1$. Consequently
\begin{equation}\label{eq:average}
                 X\in R_1'\quad\Longrightarrow\quad
                 Y=\sum_k v_kXv_k^*\in L'=P_{-1}.
\end{equation}
Let $\ell\in L$. Expanding $\ell v_k$ in this basis and moving its coefficients past $X$ gives
\[
 \ell Y=\sum_jv_jX\sum_k E_{R_1}(v_j^*\ell v_k)v_k^*
       =\sum_jv_jXv_j^*\ell=Y\ell.
\]

Let $\Omega=1\in\mathcal H$ and let $\Tr$ be the ordinary operator
trace. Choose a positive invertible $c\in Z(P_1)$ with
$\tau(T)=\Tr(cT)$ on $P_1$. Such a $c$ is obtained by rescaling the
operator trace on each matrix summand. For $x\in R_1'\cap P_1$ apply
\eqref{eq:average} to $X=cx$. The vectors $u_k^*\Omega$ form an
orthonormal basis of $\mathcal H$ since $K=\rho(P_0)$ and
$\tau^J(\rho(p))=\tau(p)$ for $p\in P_0$. Thus, for $p\in P_0$,
\[
 \tau(pY)=\sum_k
 \langle p a^{1/2}c x a^{1/2}u_k^*\Omega,u_k^*\Omega\rangle
   =\Tr(p a^{1/2}c x a^{1/2})=\tau(pax).
\]
Here $p$ commutes with $K$, and $x$ commutes with $a\in R_1$.
Thus $E_{P_0}(ax)=Y\in P_{-1}$. Taking $x=1$ and using $a\in P_0'\cap Q_1$ gives $z=E_{P_0}(a)\in P_{-1}\cap Q_0\cap Z(P_0)$. Since $a\geq\eta 1$ for some $\eta>0$, we have $z=E_{P_0}(a)\geq\eta 1$. Put $d=z^{-1}a$ and define
\[
                         F(x)=E_{P_0}(dx),\qquad x\in P_1.
\]
Then $d$ is positive and invertible in $P_0'\cap Q_1$ and $E_{P_0}(d)=1$. Testing against $P_0$ under $\tau$ gives $F(x)=E_{P_0}(d^{1/2}xd^{1/2})$, so $F$ is a faithful $P_0$-bimodular conditional expectation. Since $d\in Q_1$, we have $F(Q_1)\subset Q_0$. For $x\in R_1'\cap P_1$, we have $F(x)=z^{-1}E_{P_0}(ax)\in P_{-1}$. Finally, $F$ fixes both $Q_0$ and $P_{-1}$, giving the two equalities.
\end{proof}

\section{Isolation and finiteness for SRC commuting squares}\label{sec:finiteness}

The following proposition adapts the normalization argument of
\cite[Lemma~1.8]{N06}, recalled in \cite[Proposition~2.9]{Ni}, to the norm
$\|x\|_2=\tau(x^*x)^{1/2}$.

\begin{proposition}\label{pr:normalize}
Let $\mathcal C$ and $\mathcal C_m$ be the commuting squares \eqref{eq:square} and \eqref{eq:moving-square}, where $u_m\in\U(Q_{-1}'\cap P_0)$ and $\|u_m-1\|_2\to0$. Assume that $\mathcal C_m$ is not isomorphic to $\mathcal C$ for any $m\geq1$. After passage to a subsequence there exist unitaries
\[
\begin{array}{ll}
 p_m\in\U(Q_{-1}'\cap P_{-1}), & p'_m\in\U(P_{-1}'\cap P_0),\\
 q_m\in\U(Q_{-1}'\cap Q_0),    & q'_m\in\U(Q_0'\cap P_0)
\end{array}
\]
such that $\widehat u_m=p_mp'_mu_mq'_mq_m$ tends to $1$, $\widehat u_m\ne1$, and
\begin{equation}\label{eq:direction}
 h=\lim_{m\to\infty}\frac{\widehat u_m-1}{i\|\widehat u_m-1\|_2}
 \ \in Q_{-1}'\cap P_0,\qquad h=h^*,\qquad \|h\|_2=1 .
\end{equation}
The direction $h$ satisfies
\begin{equation}\label{eq:direction-orthogonal}
 \begin{aligned}
 E_{P_{-1}'\cap P_0}(h)&=E_{Q_{-1}'\cap P_{-1}}(h)=0,\\
 E_{Q_0'\cap P_0}(h)&=E_{Q_{-1}'\cap Q_0}(h)=0 .
 \end{aligned}
\end{equation}
Replacing $u_m$ by $\widehat u_m$ in \eqref{eq:moving-square} gives a commuting square isomorphic to $\mathcal C_m$.
\end{proposition}

\begin{proof}
Put $G_l=\U(Q_{-1}'\cap P_{-1})\,\U(P_{-1}'\cap P_0)$ and $G_r=\U(Q_0'\cap P_0)\,\U(Q_{-1}'\cap Q_0)$. The two factors in either product commute, so these are compact groups. Right multiplication by $G_r$ fixes the fourth corner, and conjugation by an element of $G_l$ preserves $P_{-1}$ and fixes $Q_{-1}$ pointwise. Choose $\widehat u_m$ in $G_lu_mG_r$ minimizing $t_m=\|\widehat u_m-1\|_2$. Then $0<t_m\leq\|u_m-1\|_2$ by the nonisomorphism assumption. Set $h_m=(\widehat u_m-1)/(it_m)$. Since $\|h_m\|_2=1$ and
$Q_{-1}'\cap P_0$ is finite dimensional, a subsequence of $(h_m)$ converges
in $2$-norm to some $h\in Q_{-1}'\cap P_0$ with $\|h\|_2=1$.
The $2$-norm and the operator norm are equivalent on $P_0$, so along this
subsequence $h_m\to h$ and $\widehat u_m\to1$ also in operator norm.
By unitarity,
\[
 h_m^*=\frac{1-\widehat u_m^*}{it_m}
      =\widehat u_m^*\frac{\widehat u_m-1}{it_m}
      =\widehat u_m^*h_m.
\]
Taking limits gives $h^*=h$, proving \eqref{eq:direction}.

Let $a=a^*$ belong to either of the two algebras defining $G_l$. For $s\in\mathbb R$, minimality gives $\|e^{-ist_ma}\widehat u_m-1\|_2\geq t_m$. Dividing by $t_m$ and passing to the limit gives $\|h-sa\|_2\geq\|h\|_2$ for every real $s$. Hence $\tau(ha)=0$. Right multiplication by $e^{-ist_ma}$ gives the same conclusion for either algebra defining $G_r$. Since $E_X$ is the orthogonal projection onto $X$ in $L^2(P_0,\tau)$, these are exactly the four equalities in \eqref{eq:direction-orthogonal}.
\end{proof}

\begin{proof}[Proof of Theorem~\ref{th:fin}\,\textup{(i)}]
Let $\mathcal C_m$ be as in \eqref{eq:moving-square}. If infinitely many $\mathcal C_m$ fail to be isomorphic to $\mathcal C$, pass to that subsequence and apply Proposition~\ref{pr:normalize}. Use the isomorphic commuting squares given by that proposition, and again denote their implementing unitaries by $u_m$. Thus $t_m=\|u_m-1\|_2>0$ and its direction $h$ satisfies \eqref{eq:direction} and \eqref{eq:direction-orthogonal}.

For a unitary $u\in\U(Q_{-1}'\cap P_0)$ write $Q_0^u=uQ_0u^*$. The diagram \eqref{eq:moving-square} with $u$ in place of $u_m$ is a commuting square precisely when
\begin{equation}\label{eq:squareparameters}
           E_{P_{-1}}(uqu^*)=E_{Q_{-1}}(q)\qquad(q\in Q_0).
\end{equation}
For each solution $u$, define a map on the spanning vectors of $L^2(P_0,\tau)$ by
\begin{equation}\label{eq:unitarylift}
             w_u(qb)=uqu^*b\qquad(q\in Q_0,\ b\in P_{-1}).
\end{equation}
By \eqref{eq:squareparameters} the proposed map preserves inner products on these spanning vectors. To see this, apply $E_{P_{-1}}$ to $uq'^*qu^*$ in the inner product of $uqu^*b$ and $uq'u^*b'$. Thus $w_u$ is an isometry on a finite-dimensional space and is onto. Its image is $Q_0^u\cdot P_{-1}$, so the resulting commuting square is nondegenerate. The unitary $w_u$ commutes with the right action of $P_{-1}$ and fixes $L^2(P_{-1},\tau)$ pointwise. Hence $w_u\in P_1$ and $w_u e=ew_u$. It implements $q\mapsto uqu^*$ on the left action of $Q_0$.

Set $w_m=w_{u_m}$. On the spanning products in \eqref{eq:unitarylift},
\[
          \frac{w_m-1}{it_m}(qb)
          =\frac{u_m q u_m^*-q}{it_m}\,b
          \longrightarrow [h,q]b .
\]
Consequently $(w_m-1)/(it_m)$ converges in operator norm to an element
$g\in P_1$, since the space is finite dimensional.
As in \eqref{eq:direction}, self-adjointness follows from unitarity. Since $u_m^*w_m$ commutes with $Q_0$, taking the normalized difference from the identity gives
\begin{equation}\label{eq:lift}
                         g-h\in Q_0'\cap P_1 .
\end{equation}

The unitary $w_m$ fixes $e$ under conjugation. Therefore
\begin{equation}\label{eq:moving-lattice}
 Q_1^{(m)}=\langle u_m Q_0 u_m^*,e\rangle=w_m Q_1w_m^*,
 \qquad R_1^{(m)}=P_{-1}'\cap Q_1^{(m)} .
\end{equation}
Up to unitary conjugacy in $P_1$ there are only finitely many unital subalgebra types, since their matrix block sizes and inclusion multiplicities are bounded. Taking a subsequence of a single type and a convergent subsequence of its conjugating unitaries, we obtain a subalgebra $R\subset P_1$ and unitaries $v_m\in P_1$ with
\[
                  R_1^{(m)}=v_m Rv_m^*,\qquad v_m\longrightarrow1.
\]
Passing to the limit gives $R\subset P_{-1}'\cap Q_1$. This inclusion suffices for the rest of the proof.

For each $m$, choose a conditional expectation $F_m:P_1\to P_0$ satisfying \eqref{eq:expectation-images} for the moving commuting square and its extension \eqref{eq:moving-lattice}. These maps have norm one, so a subsequence converges in norm to a $P_0$-bimodular map $F$ fixing $P_0$. We claim that
\begin{equation}\label{eq:limitimages}
                 F(Q_1)\subset Q_0,\qquad F(R'\cap P_1)\subset P_{-1}.
\end{equation}
For $x\in Q_1$, the elements $F_m(w_mxw_m^*)$ belong to $Q_0^{(m)}$ and converge to $F(x)$. For $x\in R'\cap P_1$, we have $v_mxv_m^*\in(R_1^{(m)})'\cap P_1$, so $F_m(v_mxv_m^*)\in P_{-1}$. This proves \eqref{eq:limitimages}.

Fix $r\in\U(R)$ and put $r_m=v_mrv_m^*$. Then $r_m\in P_{-1}'\cap P_1$ and $w_m^*r_mw_m\in Q_1$. Subtracting and dividing by $it_m$ gives
\[
 \frac{r_m-w_m^*r_mw_m}{it_m}\ \longrightarrow\ [g,r]
                 \ \in (P_{-1}'\cap P_1)+Q_1 .
\]
Indeed the numerator equals $(1-w_m^*)r_mw_m+r_m(1-w_m)$, so only $r_m\to r$ is needed in taking the limit. Both summands on the right are stable under right multiplication by $r^*$. Averaging $[g,r]r^*=g-rgr^*$ over $\U(R)$ yields
\begin{equation}\label{eq:tangent}
                      g\in(P_{-1}'\cap P_1)+Q_1+(R'\cap P_1).
\end{equation}
Apply $F$ to \eqref{eq:lift} and \eqref{eq:tangent}. Bimodularity and \eqref{eq:limitimages} yield
\[
                   h\in(P_{-1}'\cap P_0)+Q_0+P_{-1}+(Q_0'\cap P_0).
\]
Since $h$ commutes with $Q_{-1}$, averaging over $\U(Q_{-1})$ gives
\begin{equation}\label{eq:direction-sum}
 h\in(P_{-1}'\cap P_0)+(Q_{-1}'\cap P_{-1})
       +(Q_0'\cap P_0)+(Q_{-1}'\cap Q_0).
\end{equation}
By \eqref{eq:direction-orthogonal}, $h$ is orthogonal to each of the four summands. Hence $\|h\|_2^2=0$, contradicting \eqref{eq:direction}. This proves \textup{(i)}.
\end{proof}

\begin{proof}[Proof of Theorem~\ref{th:fin}\,\textup{(ii)}]
Fix the tracial inclusion data. If they are realized, choose one commuting square \eqref{eq:square} with those data. After identifying $Q_{-1}\subset P_{-1}\subset P_0$, every fourth corner is $uQ_0u^*$ for some $u\in\U(Q_{-1}'\cap P_0)$. Indeed, the inclusion multiplicities give a unitary $v\in P_0$ conjugating the two copies of $Q_0$ with the prescribed block labels. The maps $q\mapsto vqv^*$ and $q\mapsto q$ from $Q_{-1}$ into the target copy have the same multiplicities. Conjugating the first map by a unitary $w$ in that copy makes the two maps agree pointwise, so $wv$ commutes with $Q_{-1}$.

Let $(\mathcal C_m)$ be a sequence of SRC commuting squares with the fixed data, and choose representing unitaries $v_m$. The solution set of \eqref{eq:squareparameters} is closed in $\U(Q_{-1}'\cap P_0)$. Passing to a subsequence gives $v_m\to v$ in this set. Its limiting fourth corner $vQ_0v^*$ gives a nondegenerate commuting square by \eqref{eq:unitarylift}. Relative to this limiting commuting square, the fourth corners are implemented by $u_m=v_mv^*$, which tends to $1$. Part \textup{(i)} gives eventual isomorphism along the subsequence. If a fourth corner is already prescribed as the subspace limit, every such subsequence has that same limiting corner. Applying the argument to any subsequence of terms not isomorphic to the prescribed limit gives a contradiction.

It follows that a sequence of pairwise nonisomorphic SRC commuting squares with the fixed data must be finite. For a fixed tracial algebra $(P_0,\tau)$ there are only finitely many choices of block sizes and inclusion matrices for the subalgebras, and these determine their restricted traces. When $P_{-1}\subset P_0$ is connected, its inclusion matrix and block sizes determine the normalized Markov trace by the Perron--Frobenius equations \cite[Ch.~2]{GHJ}. This proves the remaining assertions.
\end{proof}

\begin{remark}
If the upper inclusion is disconnected, varying the trace can give infinitely many
pairwise nonisomorphic commuting squares even when all four algebras and their
inclusions are fixed. Take all four algebras to be $\C^2$, with identity
inclusions and trace $\tau_t(a,b)=ta+(1-t)b$ for $0<t<1/2$.
Each commuting square satisfies \eqref{eq:src}, since $R_1=Q_0=Q_1=\C^2$.
A trace-preserving isomorphism preserves the unordered pair $\{t,1-t\}$ of
traces of the minimal projections. Since $t$ is the smaller weight, distinct
values of $t$ give nonisomorphic commuting squares.
\end{remark}

\section{Standard commuting squares and Ocneanu rigidity}\label{sec:depth}

We now start with a given finite-index inclusion
$N=M_{-1}\subset M=M_0$ of $\mathrm{II}_1$ factors and its Jones tower
$(M_j)$, with $e_j\in M_j$. We use the finite-dimensional algebras in
its standard invariant to form a commuting square. When $N\subset M$ has depth
at most two, the triple
\[
          N'\cap M\ \subset\ N'\cap M_1\ \subset\ N'\cap M_2
\]
is a basic construction with Jones projection $e_2$. This includes
depth one, as in $R\otimes 1\subset R\,\overline\otimes\,M_n(\C)$
for a $\mathrm{II}_1$ factor $R$ and $n\geq2$ \cite[\S4]{BG}.
The next proposition gives the spanning identity in the standard
invariant of the given inclusion. It will let us apply the finiteness
theorem after grouping levels in any finite-depth Jones tower.

\begin{proposition}\label{pr:depth}
Let $N\subset M$ be a finite-index inclusion of $\mathrm{II}_1$ factors satisfying the depth~$2$ condition. In its Jones tower consider the standard commuting square
\begin{equation}\label{eq:standard-square}
  \mathcal C_{N,M}=
  \left(\begin{array}{ccc}
        N'\cap M_1&\subset&N'\cap M_2\\
        \cup&&\cup\\
        M'\cap M_1&\subset&M'\cap M_2
  \end{array},\,\tau\right).
\end{equation}
Then $\mathcal C_{N,M}$ has a spanning relative commutant. More explicitly,
\[
                    (M_1'\cap M_3)\cdot(M'\cap M_2)=M'\cap M_3 .
\]
\end{proposition}
\begin{proof}
We use the same Jones-projection identity at two consecutive levels.
For $i=0,1$ put $D_i=M_i'\cap M_{i+2}$. The depth-two condition
passes to the successive inclusions in the tower
\cite[\S3]{BG}, \cite[\S5]{Ni}. Hence
$M_i'\cap M_{i+3}=\spn D_i e_{i+3}D_i$.
Choose $b_a,c_a\in D_i$ with $1=\sum_a b_a e_{i+3}c_a$.
For $x\in M_{i-1}'\cap M_{i+2}$, multiply on the right by
$xe_{i+3}$ and use the Jones relation. Since $y\mapsto ye_{i+3}$
is injective on $M_{i+2}$, this gives
\[
             x=\sum_a b_a E_{M_{i+1}}(c_a x).
\]
Each coefficient commutes with $M_{i-1}$ by bimodularity. Thus
\[
      D_i\cdot(M_{i-1}'\cap M_{i+1})=M_{i-1}'\cap M_{i+2}.
\]
For $i=0$ this proves nondegeneracy of \eqref{eq:standard-square}
after taking adjoints. Both rows extend by basic construction with
the common projection $e_3$ and Markov scalar $[M:N]^{-1}$, giving
$P_1=N'\cap M_3$ and $Q_1=M'\cap M_3$.
Since $e_1\in P_{-1}=N'\cap M_1$ and $M_1=\langle M,e_1\rangle$,
we have $P_{-1}'\cap Q_1=M_1'\cap M_3=D_1$. The identity for $i=1$
is therefore SRC.
\end{proof}

\begin{proof}[Proof of Theorem~\ref{th:oc}]
Let $G$ be the adjacency matrix of the prescribed principal graph $\Gamma$, with root $*$. Let $\ell$ be the greatest distance from the root to an odd vertex. Put
\[
                   B_0=M_{\ell-1},\qquad B_1=M_{2\ell-1},
                   \qquad B_2=M_{3\ell-1}.
\]
By the multistep basic construction \cite{PP}, $B_1,B_2$ are the first two basic constructions for $N\subset B_0$.

In this convention, the size of the block of $N'\cap M_{j-1}$ indexed by a vertex $v$ at level $j$ is the number $(G^j)_{*,v}$ of paths of length $j$ from $*$ to $v$ \cite[Theorems~3.1.3 and~4.2.1]{Jo}. Every odd vertex occurs at level $\ell$ by backtracking. Every even
vertex has an odd neighbor, so its distance from the root is at most
$\ell+1\le 2\ell$ and it occurs at level $2\ell$. Thus the inclusion matrices for $N'\cap B_0\subset N'\cap B_1\subset N'\cap B_2$ are the odd-to-even block of $G^\ell$ and its transpose. To see that this triple is a basic construction, write $A_i=N'\cap B_i$ and let $f\in B_2$ be the Jones projection for $B_0\subset B_1$. It commutes with $N$. The equality $E_{B_0}(A_1)=A_0$ shows that $f$ restricts to the Jones projection from $L^2(A_1)$ onto $L^2(A_0)$. Thus restriction to $L^2(A_1)$ maps $A_1 f A_1\subset A_2$ onto the basic construction of $A_0\subset A_1$, where $A_1 f A_1$ denotes the linear span of the products. By the transpose rule, the blocks of the latter algebra have sizes $(G^{3\ell})_{*,v}$ for odd $v$. Their squared sizes sum to $\dim_\C A_2$. Consequently $A_1 f A_1=A_2$ and restriction is faithful. Thus $N\subset B_0$ satisfies the depth~$2$ condition.

Apply Proposition~\ref{pr:depth} to $N\subset B_0$. Its standard commuting square is
\begin{equation}\label{eq:grouped-square}
  \mathcal C_{N,B_0}=
  \left(\begin{array}{ccc}
       N'\cap B_1&\subset&N'\cap B_2\\
       \cup&&\cup\\
       B_0'\cap B_1&\subset&B_0'\cap B_2
  \end{array},\,\tau\right).
\end{equation}
The upper inclusion in \eqref{eq:grouped-square} is connected. Its
matrix counts length-$\ell$ walks from even to odd vertices, and every
edge of $\Gamma$ extends to such a walk with the same endpoints by
backtracking. Its upper-right dimension is
\[
               \dim_\C(N'\cap B_2)
               =\sum_{v\ {\rm odd}}(G^{3\ell})_{*,v}^{\,2}
               =(G^{6\ell})_{*,*}.
\]
Thus Theorem~\ref{th:fin}\,\textup{(ii)} leaves finitely many isomorphism classes for $\mathcal C_{N,B_0}$ in \eqref{eq:grouped-square}. Each such commuting square determines the standard invariant of $N\subset B_0$. Indeed, both rows of \eqref{eq:grouped-square} continue by horizontal basic construction with a common Jones projection at each step. Iterating \eqref{eq:grouped-square} therefore determines, up to trace-preserving isomorphism, every later commuting square in these two rows \cite[\S1.1.6]{Pa}. To identify the rows in Popa's reconstruction theorem, extend $B_0\subset B_1\subset B_2$ to the Jones tower $(B_j)_{j\geq 0}$ of $N\subset B_0$. At every level $j\geq1$ one has
\begin{equation}\label{eq:shifted-row}
 B_1'\cap B_j=(N'\cap B_1)'\cap (B_0'\cap B_j).
\end{equation}
Indeed $N'\cap B_1$ contains the Jones projection which, together with $B_0$, generates $B_1$. The right side of \eqref{eq:shifted-row} is determined inside $N'\cap B_j$ by the iterated commuting square. Hence we recover the two rows $B_0'\cap B_j$ and $B_1'\cap B_j$, with their inclusions and traces. Popa's reconstruction theorem applied to two sufficiently late consecutive columns in these rows now determines the standard invariant of $N\subset B_0$ (see \cite{Po} and \cite[\S1.3.7, p.~185]{Pa}).

It remains to recover the original inclusion. Write $\ell=2r+1$ and set $p=e_2e_4\cdots e_{2r}\in M'\cap B_0$, taking $p=1$ if $r=0$. The factors commute, and successive Markov expectations give $\tau(p)=[M:N]^{-r}>0$. Successive applications of $e_{2j}M_{2j}e_{2j}=M_{2j-2}e_{2j}$ give $pB_0p=Mp$. The normalized trace of $pB_0p$ restricts to the unique normalized trace on the factor $M$ via $x\mapsto xp$. Thus this map identifies $N\subset M$ with $Np\subset pB_0p$, with its normalized corner trace. The inclusion $N\subset B_0$ is extremal because it has finite depth \cite[\S1.3.6, p.~183]{Pa}. By Jones's reduction theorem \cite[Corollary~4.2.16 in the published version, Corollary~4.2.14 in the preprint]{Jo}, the standard invariant of $N\subset B_0$ together with $p\in N'\cap B_0$ determines the standard invariant of the corner. Unitarily conjugate projections in $N'\cap B_0$ give isomorphic corners. There are only finitely many possible block ranks for a projection in this finite-dimensional algebra, which proves the theorem.
\end{proof}

\begin{corollary}[Ocneanu's Rigidity Theorem]\label{cor:oc-rigidity}
Fix the fusion rules of a finite-depth subfactor paragroup. For each finite-index inclusion $N\subset M$ of $\mathrm{II}_1$ factors whose paragroup realizes these rules, take as generator the direct sum of one representative of each irreducible $N$--$M$ label. Only finitely many isomorphism classes of standard invariants arise in this way. Equivalently, there are only finitely many systems of unitary quantum $6j$-symbols for these rules, up to unitary changes of bases in the multiplicity spaces.
\end{corollary}

\begin{proof}
Write $I_{00}$ and $I_{01}$ for the irreducible $N$--$N$ and $N$--$M$ labels, respectively. For $a\in I_{00}$ and $i,b\in I_{01}$, let $N_{ai}^{b}$ denote the multiplicity of $b$ in the product $ai$. Choose $X=\bigoplus_{i\in I_{01}}i$. Every irreducible summand of the original generator occurs in $X$. Hence $X$ generates the whole paragroup, since the alternating products of $X$ and its dual contain, up to direct sums and summands, all the original labels.

The standard invariant with generator $X$ has even vertices $I_{00}$, odd vertices $I_{01}$, and root the $N$--$N$ identity label. The number of edges from $a$ to $b$ is the multiplicity of $b$ in $aX$, namely $\sum_{i\in I_{01}}N_{ai}^{b}$ \cite[\S2]{Kaw}. Thus every realization gives the same finite rooted principal graph. The paragroup with generator $X$ determines a finite-depth unitary planar algebra from which the pair can be recovered \cite[Corollary~2.23]{BCEGP}. We use its standard spherical structure. Since the identity objects at $N$ and $M$ are simple, both $0$-box spaces are one-dimensional, so the planar algebra is connected. Jones's realization theorem \cite[Theorem~4.3.1]{Jo} therefore gives a finite-index inclusion of $\mathrm{II}_1$ factors with this standard invariant (see also \cite[\S2.5]{BCEGP}). Theorem~\ref{th:oc} gives only finitely many possibilities for the lattice of higher relative commutants. Each lattice determines its paragroup by reconstruction. With the fusion rules fixed, choosing orthonormal bases in the multiplicity spaces expresses its associativity maps as quantum $6j$-symbols \cite[\S2]{Kaw}. Changes of these bases give precisely the unitary gauge equivalence for fixed labels. There are only finitely many ways to label the simple summands.
\end{proof}

\subsection*{The inclusion obtained by iterating the rows}
Suppose here that both row limits are also $\mathrm{II}_1$ factors.
Their inclusion is the column subfactor of the transposed commuting square, in
which $P_{-1}$ and $Q_0$ exchange roles. We call this inclusion the
\emph{row subfactor} in the following corollary.

\begin{corollary}\label{cor:src-presentation}
A subfactor obtained by iterating the rows of a commuting square in this
setting has finite depth if and only if it admits a row presentation by
a commuting square satisfying SRC.
In particular, iterating the rows of SRC commuting squares gives subfactors of
arbitrarily large finite depth.
\end{corollary}
\begin{proof}
If a commuting square satisfies SRC, its column subfactor has depth
at most two by Theorem~\ref{th:src-depth}. Sato's theorem \cite[Corollary~2.2]{SatoSquare} gives finite
depth for the row subfactor.

Conversely, start with a commuting square whose row subfactor has finite depth.
Its column subfactor also has finite depth by the same theorem,
so \eqref{eq:higher-src} gives a $j\geq1$ such that $(\mathrm{SRC}_m)$ holds for every
$m\geq j$. Choose an odd $\ell\geq j$. Retain columns
$-1,\ell-1,2\ell-1,3\ell-1,\ldots$ in the two horizontal towers.
By composition, the first two retained columns form a nondegenerate
commuting square. Its row towers are exactly the retained
towers by the multistep basic construction
\cite{PP} (see also \cite[Facts~2.20]{BCEGP}). The multistep Jones
projections are common to the two rows since they are given by the
same words in the original Jones projections. The odd step preserves
connectedness. On each horizontal inclusion graph every edge extends
to a walk of length $\ell$ with the same endpoints by backtracking,
and the vertical inclusion matrices in the grid repeat with period two.

For $m=\ell,\ldots,2\ell-1$, use $Q_m=R_m\cdot Q_{m-1}$ successively.
The algebras $R_m$ are nested and unital, so this gives
\[
       Q_{2\ell-1}=R_{2\ell-1}\cdot Q_{\ell-1}
       =(P_{-1}'\cap Q_{2\ell-1})\cdot Q_{\ell-1}.
\]
This is SRC for the commuting square formed by columns $-1$ and $\ell-1$.
The retained algebras are cofinal in each row, so their weak closures
give the original inclusion. Applying this construction to the
Temperley--Lieb--Jones subfactors with principal graphs $A_n$ rooted
at an endpoint gives SRC presentations of row depth $n-1$
\cite[\S4.7]{GHJ}.
\end{proof}

\section{Depth two for Hadamard spin models}\label{sec:spin}

For a complex Hadamard matrix $H$ of order $n>1$, Theorem~\ref{th:src-depth}
shows that the associated spin subfactor has depth two precisely when its
commuting square satisfies SRC. Here SRC reduces to $\dim_{\C} R_1=n$,
which holds exactly when $H$ is equivalent to the character table of a
finite abelian group.

Let $H=(h_{ij})$ be of order $n$, with $|h_{ij}|=1$ and $HH^*=nI$. Write $U=n^{-1/2}H$ and let $D_n$ denote the diagonal algebra. In the notation of \eqref{grid}, the spin model commuting square associated to $H$ by Popa's orthogonal-pair construction \cite{Orth} is
\begin{equation}\label{eq:spin-square}
  \mathcal C(H)=
  \left(\begin{array}{ccc}
    P_{-1}=D_n &\subset& P_0=M_n(\C)\\
    \cup && \cup\\
    Q_{-1}=\C  &\subset& Q_0=UD_nU^*
  \end{array},\,\tau\right),
\end{equation}
where $\tau$ is the normalized matrix trace. Two Hadamard matrices are \emph{equivalent} if one is obtained from the other by permuting rows and columns and multiplying individual rows and columns by scalars of modulus one. Two such matrices are equivalent if and only if their associated commuting squares are isomorphic \cite[\S2, p.~455]{NH}. Indeed, for $V=K/\sqrt n$, a unitary conjugating the two commuting squares
must normalize $D_n$ and take $UD_nU^*$ to $VD_nV^*$. Since the
normalizer of $D_n$ consists of monomial unitaries, this is equivalent
to $V=WUZ$ for monomial unitaries $W,Z$.

We dephase $H$ so that its first row and column consist of ones. Let $F_m$ denote the Fourier matrix of order $m$, with entries $\exp(2\pi i jk/m)$ for $0\le j,k<m$.

\begin{theorem}\label{th:spin}
Let $H$ be a complex Hadamard matrix of order $n\geq1$, and let $r_1,\ldots,r_n$ be the rows of a dephased representative of $H$. In the basic construction of $\mathcal C(H)$ in \eqref{eq:spin-square}, put $R_1=D_n'\cap Q_1$. The following are equivalent.
\begin{enumerate}
 \item[\textup{(i)}] $\mathcal C(H)$ satisfies \eqref{eq:src}.
 \item[\textup{(ii)}] $\dim_\C R_1=n$.
 \item[\textup{(iii)}] The entrywise quotient of any two rows is again a row, that is,
 \begin{equation}\label{eq:rowquotient}
           r_i\overline{r_j}\in\{r_1,\ldots,r_n\}
                    \qquad(1\leq i,j\leq n).
 \end{equation}
 \item[\textup{(iv)}] $H$ is equivalent to the character table of an abelian group of order $n$, or equivalently to a tensor product of Fourier matrices.
\end{enumerate}
\end{theorem}

\begin{proof}
Since $D_n$ is maximal abelian and \eqref{eq:spin-square} is a commuting square, $R_0=D_n'\cap Q_0=\C1$. We have $\dim_\C Q_0=n$, and $Q_1$, the basic construction for $\C\subset Q_0$, is isomorphic to the full matrix algebra on $L^2(Q_0,\tau)$. Thus $\dim_\C Q_1=n^2$, and \eqref{eq:size} gives \textup{(i)}$\Leftrightarrow$\textup{(ii)}.
The equivalence of \textup{(ii)} and \textup{(iv)} is the classical
maximal-dimension criterion of Jaeger, Matsumoto and Nomura
\cite[\S5.1]{JMN}, stated for spin models by Jones immediately after
\cite[Theorem~2.11.15]{Jo}.

Under \textup{(iii)} the dephased rows form an abelian group $G$ under
entrywise multiplication. Each column gives a character by coordinate
evaluation, and column orthogonality gives $n$ distinct characters.
These exhaust $\widehat G$, proving \textup{(iv)}. Conversely, dephasing
a character table at a row $g_0$ and a column $\chi_0$ gives entries
$(\chi\chi_0^{-1})(gg_0^{-1})$. Row and column phases cancel under
dephasing, so every dephased equivalent matrix has the same row group
up to relabeling. This gives \textup{(iii)} and shows that the matrix
determines $G$ up to isomorphism. Decomposing $G$ into cyclic factors
expresses its table as a tensor product of Fourier matrices.
\end{proof}

\begin{corollary}\label{cor:spin-depth}
Let $H$ be a complex Hadamard matrix of order $n>1$. The associated
spin subfactor has depth two if and only if $H$ is equivalent to the
character table of a finite abelian group. The same equivalence holds
for the inclusion obtained by iterating the rows of $\mathcal C(H)$.
\end{corollary}
\begin{proof}
The spin subfactor is irreducible of index $n$ \cite[\S2]{NH}. The first
assertion follows from Theorems~\ref{th:src-depth} and~\ref{th:spin}.
Conjugation by $U^*$ carries the transposed commuting square to $\mathcal C(H^*)$
\cite[\S2]{NH}. Taking adjoints preserves Hadamard equivalence and carries
an abelian character table to a character table of the dual group up to
relabeling, giving the second assertion.
\end{proof}

\smallskip
\noindent\textit{Real and Butson matrices.}
Write $BH(n,q)$ for the Hadamard matrices of order $n$ whose entries are $q$th roots of unity \cite{Bu}. The SRC classes in $BH(n,q)$ are represented by the character tables of abelian groups of order $n$ and exponent dividing $q$. Indeed, dephasing preserves $q$th roots of unity, so every row has order dividing $q$. Conversely, every character of such a group takes values in the $q$th roots of unity. For $q=2$ the groups are $(\mathbb Z/2\mathbb Z)^k$, so the real SRC matrices are precisely the Walsh matrices $F_2^{\otimes k}$ up to row and column permutations and sign changes.

\smallskip
\noindent\textit{Examples.}
In the classical order-four family \cite[eq.~(3) and Theorem~1]{Cr}, written in the dephased parametrization $F_4(a)$ of \cite[\S4.1]{NH}, the rows are
\[
 (1,1,1,1),\quad (1,a,-1,-a),\quad (1,-1,1,-1),\quad (1,-a,-1,a),
 \qquad |a|=1 .
\]
Squaring the entries of the second row gives $(1,a^2,1,a^2)$. For the rows to form a group, this must be the first or third row, so $a^2=\pm1$. For $a=\pm1$ the rows form the Klein four-group, and for $a=\pm i$ the second row generates all four. These are exactly the SRC parameters. For $a_m=e^{i/m}$, put $U(a)=\tfrac12F_4(a)$ and $u_m=U(a_m)U(1)^*$. Then $u_m\to1$ and $\mathcal C(F_4(a_m))$ is exactly the commuting square \eqref{eq:moving-square} with $Q_0=U(1)D_4U(1)^*$. Its limit $\mathcal C(F_4(1))$ satisfies SRC, and none of the approximating commuting squares does.
For square-free $n$, the only abelian group of order $n$ is cyclic, so $F_n$ is the unique SRC class. In particular, $BH(6,3)$ contains no SRC matrix. In order $12$ the two SRC classes are $F_{12}$ and $F_2\otimes F_2\otimes F_3$. Only the latter has a representative in $BH(12,6)$, and neither has a real representative.

\section{Quantitative estimates for SRC commuting squares}\label{sec:quantitative}

The same function will bound the number of SRC commuting squares for a given
upper-right tracial algebra and the number of standard invariants for
a given principal graph. Put
\[
                   B(1)=1,\qquad B(D)=(3D)^{25D}\quad(D\ge2).
\]

\begin{theorem}\label{th:quant-squares}
Let $D$ be a positive integer. There are at most $B(D)$ isomorphism classes of SRC commuting squares of the form \eqref{eq:square} in each of the following cases.
\begin{enumerate}
 \item[\textup{(i)}] The tracial algebra $(P_0,\tau)$ is prescribed and $\dim_\C P_0\leq D$.
 \item[\textup{(ii)}] The inclusion $P_{-1}\subset P_0$ is connected and $\dim_\C P_0\leq D$.
\end{enumerate}
In both cases the commuting squares are nondegenerate and their horizontal
inclusions have a common Markov trace. In particular, among commuting squares giving
subfactors by the column construction, at most $B(D)$ isomorphism
classes present inclusions of depth at most two when $\dim_\C P_0\leq D$.
\end{theorem}

\begin{proof}
If $D=1$ there is only the scalar commuting square. Assume $D\ge2$ and first
fix the tracial inclusion data. Section~\ref{sec:finiteness} parametrizes
the commuting squares by solutions $u\in\U(Q_{-1}'\cap P_0)$ of
\eqref{eq:squareparameters}. Fix bases $(p_i)$ of $P_{-1}$ and
$(q_j)_{j=1}^b$ of $Q_0$, and put $c=\dim_\C Q_1$ and $n=bc\le D^3$.
Define
\[
  T_u:Q_1^b\longrightarrow Q_1\oplus\!\bigoplus_{i,j} P_1,\qquad
  T_u((x_j))=\left(\sum_j x_jq_j,
                    \big([p_i,w_u x_jw_u^*]\big)_{i,j}\right),
  \qquad J(y)=(y,0).
\]
Conjugation by $w_u$ carries $Q_1$ and $Q_0$ to their moving copies.
Thus SRC says that every $(y,0)$ lies in the range of $T_u$.
Equivalently, adjoining the columns of $J$ to $T_u$ preserves its rank.
In fixed bases the entries of $w_u$ are quadratic in the real and
imaginary parts of $u$ by \eqref{eq:unitarylift}, first evaluated on
a basis of products $qb'$, $q\in Q_0$, $b'\in P_{-1}$. The entries
of $T_u$ therefore have degree at most four.

Fix a rank $0\le r\le n$. To require
$\operatorname{rank}T_u=\operatorname{rank}[T_u\ J]=r$, set all
$(r+1)$-minors of the augmented matrix equal to zero, and impose
\[
                 t\sum_{\Delta}|\Delta(T_u)|^2=1
\]
with one real variable $t$, summing over the $r$-minors of $T_u$.
For $r=0$ the empty minor is one. Together with unitarity and
\eqref{eq:squareparameters} these are real polynomial equations
in at most $2D+1$ variables, of degree at most $8n+1$
(take real and imaginary parts of complex equations).

Theorem~\ref{th:fin}\,\textup{(i)} makes the isomorphism class
locally constant on the solution set. By Milnor's bound
\cite[Theorem~2]{Mil}, the $n+1$ possible ranks together give at most
\[
           (n+1)(18D^3)^{2D+1}\le (3D)^{6D+6}
\]
classes for fixed data. We have used $8n+1\le9D^3$.

It remains to count the data. Four lists of block sizes, padded to
length $D$, give at most $(D+1)^{4D}$ choices. For an edge with
inclusion matrix $\Lambda$ and block sizes $a_i,b_j$, unitality gives
$\sum_i a_i\Lambda_{ij}=b_j$. Hence $\sum_{ij}\Lambda_{ij}\le
\sum_j b_j\le D$. Encode the matrix by listing each pair $(i,j)$
with its multiplicity and padding to length $D$. All four matrices
then give at most $(D^2+1)^{4D}$ choices. Use the same labels at
shared corners, fixing the block labels of the prescribed tracial
algebra in \textup{(i)}. The upper-right trace is prescribed in \textup{(i)}
and determined by the connected upper inclusion in \textup{(ii)}.
Consequently the number of isomorphism classes of commuting squares is at most
\begin{equation}\label{eq:square-count}
 (3D)^{6D+6}(D+1)^{4D}(D^2+1)^{4D}
       \le (3D)^{18D+6}\le B(D).
\end{equation}
\end{proof}

\smallskip
\noindent\textit{A separation estimate.}
Measure the distance between two fourth corners $A,B\subset P_0$ by
\begin{equation}\label{eq:quant-distance}
                   d(A,B)=\sup_{\|x\|_2=1}\|E_A(x)-E_B(x)\|_2.
\end{equation}
For commuting squares arising from spin models, the entries of these expectations give a direct numerical test.

\begin{proposition}\label{pr:spin-separation}
Let $H$ and $K$ be complex Hadamard matrices of order $n$ whose associated commuting squares satisfy SRC. Write $A$ and $B$ for their fourth corners in \eqref{eq:spin-square}, with the same diagonal algebra $D_n$ in the upper-left corner. If $d(A,B)<1/n$, then the two commuting squares are isomorphic.
\end{proposition}

\begin{proof}
By Theorem~\ref{th:spin}, we may label the rows by elements $g_i$
of an abelian group. For a character table, character orthogonality
gives the following formula in the usual matrix units.
\[
 \big|(E_A(e_{kl}))_{ij}\big|
    =\frac1n\,\mathbf 1_{\{g_i g_j^{-1}=g_k g_l^{-1}\}} .
\]
It persists under Hadamard equivalence. Column operations preserve
$A$, row permutations relabel the indices, and row phases multiply
the entries by scalars of modulus one. The $\sqrt n\,e_{kl}$ form
an orthonormal basis of $L^2(M_n,\tau)$, in which the matrix entries of
$E_A$ are $(E_A(e_{kl}))_{ij}$. Hence $d(A,B)<1/n$ forces the same
zero pattern for both expectations. Fix a label $o$. The pattern
recovers the group law by declaring $i\cdot j=k$ exactly when
$g_i g_o^{-1}=g_k g_j^{-1}$, using $g_i g_o^{-1}$ as the element
labeled by $i$. Thus the groups are isomorphic and their character
tables are equivalent.
\end{proof}

\smallskip
\noindent\textit{The estimate for a principal graph.}
In the proof of Theorem~\ref{th:oc}, we first choose the standard commuting square \eqref{eq:grouped-square} of $N\subset B_0$, where $B_0=M_{\ell-1}$, and then a projection in $N'\cap B_0$. Counting paths on the principal graph of $N\subset M$ bounds the number of possibilities for each choice.

\begin{theorem}\label{th:quant-graph}
Let $\Gamma$ be a finite connected rooted bipartite multigraph with at least one edge, adjacency matrix $G$ and root $*$. Let $\ell$ be the greatest distance from $*$ to an odd vertex and put
\begin{equation}\label{eq:quant-graph}
                    D_\Gamma=(G^{6\ell})_{*,*}.
\end{equation}
There are at most $B(D_\Gamma)$ isomorphism classes of standard invariants of finite-index $\mathrm{II}_1$ subfactors with principal graph $\Gamma$.
\end{theorem}

\begin{proof}
Put $D=D_\Gamma$. The commuting square \eqref{eq:grouped-square} has upper-right
dimension $D$, has connected upper inclusion, and satisfies SRC. It
determines the grouped standard invariant by Section~\ref{sec:depth}.
To recover the original invariant one also chooses a projection in
$N'\cap B_0$. If this algebra has block sizes $m_v$, it has
$\prod_v(m_v+1)\le 2^{\sum_v m_v}\le 2^D$ unitary conjugacy classes
of projections, since $N'\cap B_0\subset N'\cap B_2$. Thus
\eqref{eq:square-count} bounds the number of invariants by
\[
                  2^D(3D)^{18D+6}\le (3D)^{25D}=B(D)
                     \qquad(D\ge2).
\]
If $D=1$, the only connected rooted graph that can occur is the single
edge with multiplicity one, and its index is one. For irreducible
$N\subset M$, the recovering projection is minimal since its corner
in $N'\cap B_0$ is isomorphic to $N'\cap M$. In that case one need
only count the odd vertices. For $A_3$ rooted at an endpoint one has
$\ell=1$, $D_\Gamma=4$, and just one possible corner.
\end{proof}

The estimate also gives a bound in terms of index and depth alone. Suppose $[M:N]\le I$ and the depth is at most a positive integer $h$. At finite depth one has $\|G\|^2=[M:N]$ \cite[\S1.3.6, p.~183]{Pa}. Since $\ell\le h$, we obtain $D_\Gamma\le\|G\|^{6\ell}\le I^{3h}$. One can count the possible graphs themselves as well. Since $\sum_w G_{vw}\le\sum_w G_{vw}^2\le I$, there are at most $V=1+\lfloor I\rfloor+\cdots+\lfloor I\rfloor^h$ vertices. Fix the root label and pad the adjacency matrix with isolated vertices. Each entry above the diagonal is at most $\lfloor\sqrt I\rfloor$, so there are at most $(1+\lfloor\sqrt I\rfloor)^{\binom V2}$ rooted principal graphs. Multiplying this last number by $B(\lceil I^{3h}\rceil)$ bounds the total number of standard invariants with these index and depth bounds.

\section*{Acknowledgment}
The second author is grateful to Sorin Popa for bringing this problem to his
attention and for their earlier conversations on the subject.

\section*{Funding}
This research received no specific grant from any funding agency in the public, commercial, or not-for-profit sectors.

\section*{Declaration of competing interest}
The authors declare that they have no competing interests.

\section*{Data availability}
No datasets were generated or analyzed in this theoretical work.

\section*{Declaration of generative AI and AI-assisted technologies in the manuscript preparation process}
We used OpenAI’s ChatGPT to help revise parts of the exposition. Its use in developing the arguments is described in the introduction. We take full responsibility for the mathematical results, the references, and their presentation.

\end{document}